\documentclass[12pt]{article}
\usepackage{times,mathptmx}
\usepackage{pstricks,pst-node,graphics}
\input epsf

\title{The many faces of alternating-sign matrices}

\author{James Propp\thanks{Supported by grants from the National Science
Foundation and the National Security Agency.}\\
\small Department of Mathematics\\[-0.8ex]
\small University of Wisconsin -- Madison, Wisconsin, USA\\[-0.8ex]
\small \texttt{propp@math.wisc.edu}}

\begin{document}
\maketitle

\begin{abstract}
\noindent
I give a survey of different combinatorial forms of 
alternating-sign matrices, starting with the original form introduced 
by Mills, Robbins and Rumsey as well as corner-sum matrices, 
height-function matrices, three-colorings, monotone triangles, 
tetrahedral order ideals, square ice, gasket-and-basket tilings
and full packings of loops.
(This article has been published in a conference edition
of the journal {\it Discrete Mathematics and Theoretical 
Computer Science}, entitled ``Discrete Models: Combinatorics,
Computation, and Geometry,'' edited by R.\ Cori, J.\ Mazoyer,
M.\ Morvan, and R.\ Mosseri, and published in July 2001
in cooperation with le Maison de l'Informatique et des
Math\'ematiques Discr\`etes, Paris, France: ISSN 1365-8050,
{\tt http://dmtcs.lori.fr}.)
\end{abstract}

\newcommand{\R}{{\bf R}}
\newcommand{\Z}{{\bf Z}}
\newcommand{\odisk}{\pscircle[fillstyle=solid,fillcolor=black]}

\section{Introduction}
\label{sec:in}              
An \emph{alternating-sign matrix of order $n$}
is an $n$-by-$n$ array of 0's, $+1$'s and $-1$'s
with the property that in each row and each column,
the non-zero entries alternate in sign,
beginning and ending with a $+1$.
For example, Figure 1 shows an alternating-sign matrix
(ASM for short)
of order 4.
$$
\left( \begin{array}{rrrr}
 0 & +1 &  0 &  0 \\
+1 & -1 & +1 &  0 \\
 0 &  0 &  0 & +1 \\
 0 & +1 &  0 &  0 
\end{array} \right)
$$
\begin{center}
Figure 1: An alternating-sign matrix of order 4.
\end{center}
Figure 2 exhibits all seven of the ASMs of order 3.
$$
\left( \begin{array}{rrr}
 0 &  0 & +1 \\
 0 & +1 &  0 \\
+1 &  0 &  0 
\end{array} \right)
\ \ \left( \begin{array}{rrr}
 0 &  0 & +1 \\
+1 &  0 &  0 \\
 0 & +1 &  0 
\end{array} \right)
\ \ \left( \begin{array}{rrr}
 0 & +1 &  0 \\
 0 &  0 & +1 \\
+1 &  0 &  0 
\end{array} \right)
\ \ \left( \begin{array}{rrr}
 0 & +1 &  0 \\
+1 & -1 & +1 \\
 0 & +1 &  0 
\end{array} \right)
$$
$$
\left( \begin{array}{rrr}
 0 & +1 &  0 \\
+1 &  0 &  0 \\
 0 &  0 & +1 
\end{array} \right)
\ \ \left( \begin{array}{rrr}
+1 &  0 &  0 \\
 0 &  0 & +1 \\
 0 & +1 &  0 
\end{array} \right)
\ \ \left( \begin{array}{rrr}
+1 &  0 &  0 \\
 0 & +1 &  0 \\
 0 &  0 & +1 
\end{array} \right)
$$
\begin{center}
Figure 2: The seven alternating-sign matrices of order 3.
\end{center}
Matrices satisfying these constraints
were first investigated by Mills, Robbins and Rumsey~\cite{MRR}.
The matrices arose from their investigation
of Dodgson's scheme for computing determinants via ``condensation''
(see section~\ref{sec:dodgson}).
The number of ASMs of order $n$, for small values of $n$,
goes like $1,2,7,42,429,7436,\dots\/$,
and it was conjectured by Mills et al.\ that
the number of ASMs of order $n$
is given by the product
$$\frac{1!4!7!\cdots(3n-2)!}{n!(n+1)!(n+2)!\cdots(2n-1)!}.$$
However, it took over a decade
before this conjecture was proved
by Zeilberger~\cite{Z}.
For more details on this history,
see the expository article by Robbins~\cite{Ro1},
the survey article by Bressoud and Propp~\cite{BP},
or the book by Bressoud~\cite{Br}.

Here my concern will be not with the alternating-sign matrix conjecture
and its proof by Zeilberger,
but with the inherent interest of alternating-sign matrices
as combinatorial objects
admitting many different representations. 
I will present here a number of different ways of looking at an ASM.
Along the way, I will also mention a few topics related to ASMs
in their various guises, such as 
weighted enumeration formulas
and asymptotic shape.
Much of what is in this article
has appeared elsewhere, but
I hope that by gathering these topics together in one place,
I will help raise the level of knowledge and interest
of the mathematical community
concerning these fascinating combinatorial objects.

\section{Corner-sum, heights, and colorings}
\label{sec:corner}              
Given an ASM $(a_{i,j})_{i,j=1}^{n}$ of order $n$,
we can define a \emph{corner-sum matrix} $(c_{i,j})_{i,j=0}^n$
of order $n$
by putting $c_{i,j} = \sum_{i'\leq i,\ j'\leq j} \, a_{i',j'}$.
This definition was introduced in~\cite{RR}.
Figure 3 shows the seven corner-sum matrices of order 3
(note that they are 4-by-4 matrices).
$$
\left( \begin{array}{rrrr}
 0 & 0 & 0 & 0 \\
 0 & 0 & 0 & 1 \\
 0 & 0 & 1 & 2 \\
 0 & 1 & 2 & 3
\end{array} \right)
\ \ \left( \begin{array}{rrrr}
 0 & 0 & 0 & 0 \\
 0 & 0 & 0 & 1 \\
 0 & 1 & 1 & 2 \\
 0 & 1 & 2 & 3
\end{array} \right)
\ \ \left( \begin{array}{rrrr}
 0 & 0 & 0 & 0 \\
 0 & 0 & 1 & 1 \\
 0 & 0 & 1 & 2 \\
 0 & 1 & 2 & 3
\end{array} \right)
\ \ \left( \begin{array}{rrrr}
 0 & 0 & 0 & 0 \\
 0 & 0 & 1 & 1 \\
 0 & 1 & 1 & 2 \\
 0 & 1 & 2 & 3
\end{array} \right)
$$
$$
\left( \begin{array}{rrrr}
 0 & 0 & 0 & 0 \\
 0 & 0 & 1 & 1 \\
 0 & 1 & 2 & 2 \\
 0 & 1 & 2 & 3
\end{array} \right)
\ \ \left( \begin{array}{rrrr}
 0 & 0 & 0 & 0 \\
 0 & 1 & 1 & 1 \\
 0 & 1 & 1 & 2 \\
 0 & 1 & 2 & 3
\end{array} \right)
\ \ \left( \begin{array}{rrrr}
 0 & 0 & 0 & 0 \\
 0 & 1 & 1 & 1 \\
 0 & 1 & 2 & 2 \\
 0 & 1 & 2 & 3
\end{array} \right)
$$
\begin{center}
Figure 3: The seven corner-sum matrices of order 3.
\end{center}

Corner-sum matrices, viewed as objects in their own right,
have a very simple description:
the first row and first column consist of 0's,
the last row and last column consist of the numbers from 0 to $n$
(in order),
and within each row and column,
each entry is either equal to, or one more than,
the preceding entry.

Note that the seven corner-sum matrices in Figure 3
correspond respectively to the seven alternating-sign matrices in Figure 2.
I will adhere to this pattern throughout,
to make it easier for the reader to verify
the bijections between the different representations.

Corner-sum matrices can in turn be transformed
into a somewhat more symmetrical form.
Given a corner-sum matrix $(c_{i,j})_{i,j=0}^n$
define $h_{i,j} = i+j-2c_{i,j}$.
Call the result a \emph{height-function matrix}
(see~\cite{EKLP}).
Figure 4 shows the seven height-function matrices of order 3
(4-by-4 matrices).
$$
\left( \begin{array}{rrrr}
 0 & 1 & 2 & 3 \\
 1 & 2 & 3 & 2 \\
 2 & 3 & 2 & 1 \\
 3 & 2 & 1 & 0
\end{array} \right)
\ \ \left( \begin{array}{rrrr}
 0 & 1 & 2 & 3 \\
 1 & 2 & 3 & 2 \\
 2 & 1 & 2 & 1 \\
 3 & 2 & 1 & 0
\end{array} \right)
\ \ \left( \begin{array}{rrrr}
 0 & 1 & 2 & 3 \\
 1 & 2 & 1 & 2 \\
 2 & 3 & 2 & 1 \\
 3 & 2 & 1 & 0
\end{array} \right)
\ \ \left( \begin{array}{rrrr}
 0 & 1 & 2 & 3 \\
 1 & 2 & 1 & 2 \\
 2 & 1 & 2 & 1 \\
 3 & 2 & 1 & 0
\end{array} \right)
$$
$$
\left( \begin{array}{rrrr}
 0 & 1 & 2 & 3 \\
 1 & 2 & 1 & 2 \\
 2 & 1 & 0 & 1 \\
 3 & 2 & 1 & 0
\end{array} \right)
\ \ \left( \begin{array}{rrrr}
 0 & 1 & 2 & 3 \\
 1 & 0 & 1 & 2 \\
 2 & 1 & 2 & 1 \\
 3 & 2 & 1 & 0
\end{array} \right)
\ \ \left( \begin{array}{rrrr}
 0 & 1 & 2 & 3 \\
 1 & 0 & 1 & 2 \\
 2 & 1 & 0 & 1 \\
 3 & 2 & 1 & 0
\end{array} \right)
$$
\begin{center}
Figure 4: The seven height-function matrices of order 3.
\end{center}

Height-function matrices
have a simple intrinsic description:
the first row and first column consist of the numbers from 0 to $n$
(consecutively),
the last row and last column consist of the numbers from $n$ to 0
(consecutively),
and any two entries that are row-adjacent or column-adjacent
differ by 1.

If one reduces a height-function matrix modulo 3,
and views the residues 0, 1, and 2 as ``colors'',
one obtains a proper 3-coloring of the $n+1$-by-$n+1$ square grid
satisfying specific boundary conditions.
Here ``proper'' means that adjacent sites get distinct colors,
and the specific boundary conditions are as follows:
colors increase modulo 3 along the first row and first column
and decrease modulo 3 along the last row and last column,
with the color 0 occurring in the upper left.
Figure 5 shows the seven such 3-colorings
of the 4-by-4 grid.

$$
\left( \begin{array}{rrrr}
 0 & 1 & 2 & 0 \\
 1 & 2 & 0 & 2 \\
 2 & 0 & 2 & 1 \\
 0 & 2 & 1 & 0
\end{array} \right)
\ \ \left( \begin{array}{rrrr}
 0 & 1 & 2 & 0 \\
 1 & 2 & 0 & 2 \\
 2 & 1 & 2 & 1 \\
 0 & 2 & 1 & 0
\end{array} \right)
\ \ \left( \begin{array}{rrrr}
 0 & 1 & 2 & 0 \\
 1 & 2 & 1 & 2 \\
 2 & 0 & 2 & 1 \\
 0 & 2 & 1 & 0
\end{array} \right)
\ \ \left( \begin{array}{rrrr}
 0 & 1 & 2 & 0 \\
 1 & 2 & 1 & 2 \\
 2 & 1 & 2 & 1 \\
 0 & 2 & 1 & 0
\end{array} \right)
$$
$$
\left( \begin{array}{rrrr}
 0 & 1 & 2 & 0 \\
 1 & 2 & 1 & 2 \\
 2 & 1 & 0 & 1 \\
 0 & 2 & 1 & 0
\end{array} \right)
\ \ \left( \begin{array}{rrrr}
 0 & 1 & 2 & 0 \\
 1 & 0 & 1 & 2 \\
 2 & 1 & 2 & 1 \\
 0 & 2 & 1 & 0
\end{array} \right)
\ \ \left( \begin{array}{rrrr}
 0 & 1 & 2 & 0 \\
 1 & 0 & 1 & 2 \\
 2 & 1 & 0 & 1 \\
 0 & 2 & 1 & 0
\end{array} \right)
$$
\begin{center}
Figure 5: The seven colorings associated with the ASMs of order 3.
\end{center}

Conversely, every proper 3-coloring of that graph
that satisfies the boundary conditions
is associated with a unique height-function matrix~\cite{Ba}.

\section{Monotone triangles and order ideals}
\label{sec:monotone}              
Another way to ``process'' an ASM
is to form partial sums of its columns from the top toward the bottom,
as shown in Figure 6 for a 4-by-4 ASM.
In the resulting square matrix of partial sums,
the $i$th row has $i$ 1's in it
and $n-i$ 0's.
Hence we may form a triangular array
whose $i$th row consists of precisely those values $j$
for which the $i,j$th entry of the partial-sum matrix is 1.
The result is called a \emph{monotone triangle}~\cite{MRR}
(or \emph{Gog triangle} in the terminology of Zeilberger~\cite{Z}).
Figure 7 shows the seven monotone triangles of order 3.

$$
\left( \begin{array}{rrrr}
 0 & +1 &  0 &  0 \\
+1 & -1 & +1 &  0 \\
 0 &  0 &  0 & +1 \\
 0 & +1 &  0 &  0 
\end{array} \right)
\rightarrow
\left( \begin{array}{rrrr}
0 & 1 & 0 & 0 \\
1 & 0 & 1 & 0 \\
1 & 0 & 1 & 1 \\
1 & 1 & 1 & 1 
\end{array} \right)
\rightarrow
\begin{array}{rrrrrrr}
  &   &   & 2 &   &   &   \\[2ex]
  &   & 1 &   & 3 &   &   \\[2ex]
  & 1 &   & 3 &   & 4 &   \\[2ex]
1 &   & 2 &   & 3 &   & 4
\end{array}
$$
\begin{center}
Figure 6: Turning an ASM into a monotone triangle.
\end{center}

$$
\begin{array}{rrrrr}
  &   & 3 &   &   \\[2ex]
  & 2 &   & 3 &   \\[2ex]
1 &   & 2 &   & 3
\end{array}
\ \ \ \ \ \ \ \ \begin{array}{rrrrr}
  &   & 3 &   &   \\[2ex]
  & 1 &   & 3 &   \\[2ex]
1 &   & 2 &   & 3
\end{array}
\ \ \ \ \ \ \ \ \begin{array}{rrrrr}
  &   & 2 &   &   \\[2ex]
  & 2 &   & 3 &   \\[2ex]
1 &   & 2 &   & 3
\end{array}
\ \ \ \ \ \ \ \ \begin{array}{rrrrr}
  &   & 2 &   &   \\[2ex]
  & 1 &   & 3 &   \\[2ex]
1 &   & 2 &   & 3
\end{array}
$$
$$
\begin{array}{rrrrr}
  &   & 2 &   &   \\[2ex]
  & 1 &   & 2 &   \\[2ex]
1 &   & 2 &   & 3
\end{array}
\ \ \ \ \ \ \ \ \begin{array}{rrrrr}
  &   & 1 &   &   \\[2ex]
  & 1 &   & 3 &   \\[2ex]
1 &   & 2 &   & 3
\end{array}
\ \ \ \ \ \ \ \ \begin{array}{rrrrr}
  &   & 1 &   &   \\[2ex]
  & 1 &   & 2 &   \\[2ex]
1 &   & 2 &   & 3
\end{array}
$$
\begin{center}
Figure 7: The seven monotone triangles of order 3.
\end{center}

One may intrinsically describe a monotone triangle of order $n$
as a triangular array with $n$ numbers along each side,
where the numbers in the bottom row are $1$ through $n$ in succession,
the numbers in each row are strictly increasing from left to right,
and the numbers along diagonals
are weakly increasing from left to right.
Zeilberger's proof of the ASM conjecture~\cite{Z}
used these Gog triangles
and a natural generalization,
``Gog trapezoids''.

A different geometry comes from looking at the set of ASMs
as a distributive lattice.
Given two height-function matrices $(h_{i,j})_{i,j=0}^n$
and $(h'_{i,j})_{i,j=0}^n$,
we can define new matrices
(called the \emph{join} and \emph{meet})
whose $i,j$th entries are
$\max(h_{i,j},h'_{i,j})$
and $\min(h_{i,j},h'_{i,j})$,
respectively.
These new matrices are themselves height-function matrices,
and the operations of join and meet
turn the set of ASMs of order $n$ into a distributive lattice $L$
(see~\cite{St} for background on finite posets and lattices).

The fundamental theorem of finite distributive lattices
tells us that $L$
can be realized as the lattice of order-ideals of
a certain poset $P$,
namely, the poset of join-irreducibles of the lattice $L$.
There is a nice geometric description of the ranked poset $P$.
It has $(1)(n-1)$ elements of rank 0,
$(2)(n-2)$ elements of rank 1,
$(3)(n-3)$ elements of rank 2,
etc.,
up through
$(n-1)(1)$ elements of rank $n-1$.
These elements are arranged in the fashion of a tetrahedron
resting on its edge.
A generic element of $P$, well inside the interior of the tetrahedron,
covers 4 elements and is covered by 4 elements.

Using this poset $P$,
we can give a picture of the lattice $L$
that does not require a knowledge of poset-theory
(also described in~\cite{EKLP}).
Picture a tetrahedron that is densely packed with
$(1)(n-1)+(2)(n-2)+(3)(n-3)+\dots+(n-1)(1)$ balls,
resting on an edge.
Carefully remove the two upper faces of the tetrahedron
so as not disturb the balls.
One may now start to remove some of the balls,
starting from the top,
so that removal of a ball does not affect
any of the balls below.
There are many configurations of this kind,
ranging from the full packing
to the empty packing.
These configurations are in bijective correspondence
with the ASMs of order $n$,
and the lattice operations of meet and join
correspond to intersection and union.

\section{Square ice}
\label{sec:ice}              
Zeilberger's proof of the ASM conjecture was followed in short order
by a simpler proof due to Kuperberg~\cite{Ku1}.
Kuperberg's proof made use of a different representation of ASMs,
the ``6-vertex model'' of statistical mechanics.
This model is also called square ice
on account of its origin as a two-dimensional surrogate
for a more realistic (and still intractable) 
three-dimensional model of ice proposed by physicists~\cite{Ba}.
A square ice state is an orientation of the edges of a square grid
or a finite sub-graph thereof
with the property that each vertex other than vertices on the boundary
has two incoming arrows and two outgoing arrows.
Each internal vertex must be of one of the six kinds
shown in Figure 8
(hence the name ``six-vertex model'').
The markings under the six vertex-types
can be ignored for the time being.
$$
\begin{array}{c}
\begin{pspicture}(-.9,-.9)(.9,.9)
\psset{arrowsize=5pt}
\psline{->}(0,.5)
\psline{->}(0,-.5)
\psline{-<}(.5,0)
\psline{-<}(-.5,0)
\psline(0,.9)
\psline(0,-.9)
\psline(.9,0)
\psline(-.9,0)
\end{pspicture} \\
+1
\end{array}
\begin{array}{c}
\begin{pspicture}(-1,-1)(1,1)
\psset{arrowsize=5pt}
\psline{-<}(0,.5)
\psline{-<}(0,-.5)
\psline{->}(.5,0)
\psline{->}(-.5,0)
\psline(0,.9)
\psline(0,-.9)
\psline(.9,0)
\psline(-.9,0)
\psset{arrowsize=5pt}
\end{pspicture} \\
-1
\end{array}
\begin{array}{c}
\begin{pspicture}(-1,-1)(1,1)
\psset{arrowsize=5pt}
\psline{->}(0,.5)
\psline{-<}(0,-.5)
\psline{->}(.5,0)
\psline{-<}(-.5,0)
\psline(0,.9)
\psline(0,-.9)
\psline(.9,0)
\psline(-.9,0)
\psset{arrowsize=5pt}
\end{pspicture} \\
0
\end{array}
\begin{array}{c}
\begin{pspicture}(-1,-1)(1,1)
\psset{arrowsize=5pt}
\psline{-<}(0,.5)
\psline{->}(0,-.5)
\psline{-<}(.5,0)
\psline{->}(-.5,0)
\psline(0,.9)
\psline(0,-.9)
\psline(.9,0)
\psline(-.9,0)
\psset{arrowsize=5pt}
\end{pspicture} \\
0
\end{array}
\begin{array}{c}
\begin{pspicture}(-1,-1)(1,1)
\psset{arrowsize=5pt}
\psline{->}(0,.5)
\psline{-<}(0,-.5)
\psline{-<}(.5,0)
\psline{->}(-.5,0)
\psline(0,.9)
\psline(0,-.9)
\psline(.9,0)
\psline(-.9,0)
\psset{arrowsize=5pt}
\end{pspicture} \\
0
\end{array}
\begin{array}{c}
\begin{pspicture}(-1,-1)(1,1)
\psset{arrowsize=5pt}
\psline{-<}(0,.5)
\psline{->}(0,-.5)
\psline{->}(.5,0)
\psline{-<}(-.5,0)
\psline(0,.9)
\psline(0,-.9)
\psline(.9,0)
\psline(-.9,0)
\psset{arrowsize=5pt}
\end{pspicture} \\
0
\end{array}
$$
\begin{center} 
Figure 8: The six vertex-types for the square-ice model.
\end{center}

As our finite subgraph of the square grid,
we will take the ``generalized tic-tac-toe graph''
formed by $n$ horizontal lines and $n$ vertical lines
meeting in $n^2$ intersections of degree 4,
with $4n$ vertices of degree 1 at the boundary.
We say that an ice state on this graph
satisfies \emph{domain-wall boundary conditions}~\cite{KBI}
if all the arrows along the left and right flank point inward
and all the arrows along the top and bottom point outward.

Figure 9 shows the possibilities when $n=3$.
$$
\begin{pspicture}(0,0)(2.8,2.8)
\psset{xunit=.7,yunit=.7}
\psset{arrowsize=5pt}
\newcommand{\iceR}{\psline{->}(0,0)(.7,0) \psline(1,0)}
\newcommand{\iceL}{\psline{->}(0,0)(-.7,0) \psline(-1,0)}
\newcommand{\iceU}{\psline{->}(0,0)(0,.7) \psline(0,1)}
\newcommand{\iceD}{\psline{->}(0,0)(0,-.7) \psline(0,-1)}
\rput(0,1){\iceR}
\rput(0,2){\iceR}
\rput(0,3){\iceR}
\rput(4,1){\iceL}
\rput(4,2){\iceL}
\rput(4,3){\iceL}
\rput(1,3){\iceU}
\rput(2,3){\iceU}
\rput(3,3){\iceU}
\rput(1,1){\iceD}
\rput(2,1){\iceD}
\rput(3,1){\iceD}
\rput(1,3){\iceR}
\rput(1,2){\iceR}
\rput(2,3){\iceR}
\rput(3,2){\iceL}
\rput(2,1){\iceL}
\rput(3,1){\iceL}
\rput(1,2){\iceU}
\rput(1,1){\iceU}
\rput(2,2){\iceU}
\rput(2,2){\iceD}
\rput(3,3){\iceD}
\rput(3,2){\iceD}
\end{pspicture}
\ \ \ \ \ \ \ \ {}
\begin{pspicture}(0,0)(2.8,2.8)
\psset{xunit=.7,yunit=.7}
\psset{arrowsize=5pt}
\newcommand{\iceR}{\psline{->}(0,0)(.7,0) \psline(1,0)}
\newcommand{\iceL}{\psline{->}(0,0)(-.7,0) \psline(-1,0)}
\newcommand{\iceU}{\psline{->}(0,0)(0,.7) \psline(0,1)}
\newcommand{\iceD}{\psline{->}(0,0)(0,-.7) \psline(0,-1)}
\rput(0,1){\iceR}
\rput(0,2){\iceR}
\rput(0,3){\iceR}
\rput(4,1){\iceL}
\rput(4,2){\iceL}
\rput(4,3){\iceL}
\rput(1,3){\iceU}
\rput(2,3){\iceU}
\rput(3,3){\iceU}
\rput(1,1){\iceD}
\rput(2,1){\iceD}
\rput(3,1){\iceD}
\rput(1,3){\iceR}
\rput(2,3){\iceR}
\rput(1,1){\iceR}
\rput(3,1){\iceL}
\rput(2,2){\iceL}
\rput(3,2){\iceL}
\rput(1,2){\iceU}
\rput(2,2){\iceU}
\rput(2,1){\iceU}
\rput(1,2){\iceD}
\rput(3,3){\iceD}
\rput(3,2){\iceD}
\end{pspicture}
\ \ \ \ \ \ \ \ {}
\begin{pspicture}(0,0)(2.8,2.8)
\psset{xunit=.7,yunit=.7}
\psset{arrowsize=5pt}
\newcommand{\iceR}{\psline{->}(0,0)(.7,0) \psline(1,0)}
\newcommand{\iceL}{\psline{->}(0,0)(-.7,0) \psline(-1,0)}
\newcommand{\iceU}{\psline{->}(0,0)(0,.7) \psline(0,1)}
\newcommand{\iceD}{\psline{->}(0,0)(0,-.7) \psline(0,-1)}
\rput(0,1){\iceR}
\rput(0,2){\iceR}
\rput(0,3){\iceR}
\rput(4,1){\iceL}
\rput(4,2){\iceL}
\rput(4,3){\iceL}
\rput(1,3){\iceU}
\rput(2,3){\iceU}
\rput(3,3){\iceU}
\rput(1,1){\iceD}
\rput(2,1){\iceD}
\rput(3,1){\iceD}
\rput(1,2){\iceR}
\rput(1,3){\iceR}
\rput(2,2){\iceR}
\rput(3,1){\iceL}
\rput(3,3){\iceL}
\rput(2,1){\iceL}
\rput(1,1){\iceU}
\rput(1,2){\iceU}
\rput(3,1){\iceU}
\rput(3,2){\iceU}
\rput(2,3){\iceD}
\rput(2,2){\iceD}
\end{pspicture}
\ \ \ \ \ \ \ \ {}
\begin{pspicture}(0,0)(2.8,2.8)
\psset{xunit=.7,yunit=.7}
\psset{arrowsize=5pt}
\newcommand{\iceR}{\psline{->}(0,0)(.7,0) \psline(1,0)}
\newcommand{\iceL}{\psline{->}(0,0)(-.7,0) \psline(-1,0)}
\newcommand{\iceU}{\psline{->}(0,0)(0,.7) \psline(0,1)}
\newcommand{\iceD}{\psline{->}(0,0)(0,-.7) \psline(0,-1)}
\rput(0,1){\iceR}
\rput(0,2){\iceR}
\rput(0,3){\iceR}
\rput(4,1){\iceL}
\rput(4,2){\iceL}
\rput(4,3){\iceL}
\rput(1,3){\iceU}
\rput(2,3){\iceU}
\rput(3,3){\iceU}
\rput(1,1){\iceD}
\rput(2,1){\iceD}
\rput(3,1){\iceD}
\rput(1,3){\iceR}
\rput(2,2){\iceR}
\rput(1,1){\iceR}
\rput(3,3){\iceL}
\rput(2,2){\iceL}
\rput(3,1){\iceL}
\rput(1,2){\iceU}
\rput(2,1){\iceU}
\rput(3,2){\iceU}
\rput(1,2){\iceD}
\rput(2,3){\iceD}
\rput(3,2){\iceD}
\end{pspicture}
$$
$$
\begin{pspicture}(0,0)(2.8,2.8)
\psset{xunit=.7,yunit=.7}
\psset{arrowsize=5pt}
\newcommand{\iceR}{\psline{->}(0,0)(.7,0) \psline(1,0)}
\newcommand{\iceL}{\psline{->}(0,0)(-.7,0) \psline(-1,0)}
\newcommand{\iceU}{\psline{->}(0,0)(0,.7) \psline(0,1)}
\newcommand{\iceD}{\psline{->}(0,0)(0,-.7) \psline(0,-1)}
\rput(0,1){\iceR}
\rput(0,2){\iceR}
\rput(0,3){\iceR}
\rput(4,1){\iceL}
\rput(4,2){\iceL}
\rput(4,3){\iceL}
\rput(1,3){\iceU}
\rput(2,3){\iceU}
\rput(3,3){\iceU}
\rput(1,1){\iceD}
\rput(2,1){\iceD}
\rput(3,1){\iceD}
\rput(1,1){\iceR}
\rput(2,1){\iceR}
\rput(1,3){\iceR}
\rput(3,3){\iceL}
\rput(2,2){\iceL}
\rput(3,2){\iceL}
\rput(1,2){\iceU}
\rput(3,2){\iceU}
\rput(3,1){\iceU}
\rput(1,2){\iceD}
\rput(2,2){\iceD}
\rput(2,3){\iceD}
\end{pspicture}
\ \ \ \ \ \ \ \ {}
\begin{pspicture}(0,0)(2.8,2.8)
\psset{xunit=.7,yunit=.7}
\psset{arrowsize=5pt}
\newcommand{\iceR}{\psline{->}(0,0)(.7,0) \psline(1,0)}
\newcommand{\iceL}{\psline{->}(0,0)(-.7,0) \psline(-1,0)}
\newcommand{\iceU}{\psline{->}(0,0)(0,.7) \psline(0,1)}
\newcommand{\iceD}{\psline{->}(0,0)(0,-.7) \psline(0,-1)}
\rput(0,1){\iceR}
\rput(0,2){\iceR}
\rput(0,3){\iceR}
\rput(4,1){\iceL}
\rput(4,2){\iceL}
\rput(4,3){\iceL}
\rput(1,3){\iceU}
\rput(2,3){\iceU}
\rput(3,3){\iceU}
\rput(1,1){\iceD}
\rput(2,1){\iceD}
\rput(3,1){\iceD}
\rput(1,1){\iceR}
\rput(1,2){\iceR}
\rput(2,2){\iceR}
\rput(3,3){\iceL}
\rput(2,3){\iceL}
\rput(3,1){\iceL}
\rput(2,2){\iceU}
\rput(2,1){\iceU}
\rput(3,2){\iceU}
\rput(1,2){\iceD}
\rput(1,3){\iceD}
\rput(3,2){\iceD}
\end{pspicture}
\ \ \ \ \ \ \ \ {}
\begin{pspicture}(0,0)(2.8,2.8)
\psset{xunit=.7,yunit=.7}
\psset{arrowsize=5pt}
\newcommand{\iceR}{\psline{->}(0,0)(.7,0) \psline(1,0)}
\newcommand{\iceL}{\psline{->}(0,0)(-.7,0) \psline(-1,0)}
\newcommand{\iceU}{\psline{->}(0,0)(0,.7) \psline(0,1)}
\newcommand{\iceD}{\psline{->}(0,0)(0,-.7) \psline(0,-1)}
\rput(0,1){\iceR}
\rput(0,2){\iceR}
\rput(0,3){\iceR}
\rput(4,1){\iceL}
\rput(4,2){\iceL}
\rput(4,3){\iceL}
\rput(1,3){\iceU}
\rput(2,3){\iceU}
\rput(3,3){\iceU}
\rput(1,1){\iceD}
\rput(2,1){\iceD}
\rput(3,1){\iceD}
\rput(1,1){\iceR}
\rput(1,2){\iceR}
\rput(2,1){\iceR}
\rput(3,3){\iceL}
\rput(2,3){\iceL}
\rput(3,2){\iceL}
\rput(2,2){\iceU}
\rput(3,1){\iceU}
\rput(3,2){\iceU}
\rput(1,2){\iceD}
\rput(1,3){\iceD}
\rput(2,2){\iceD}
\end{pspicture}
$$
\begin{center} 
Figure 9: The seven square-ice states for $n=3$ with domain-wall
boundary conditions.
\end{center}

These states of the square-ice model
are in bijective correspondence with ASMs.
To turn a state of the square-ice model on an $n$-by-$n$ grid
with domain-wall boundary conditions
into an alternating-sign matrix of order $n$,
replace each vertex by $+1$, $-1$, or 0
according to the marking given in Figure 8.
Kuperberg was able to give a simplified proof of the ASM conjecture
by making use of results about the square-ice model
in the mathematical physics literature~\cite{KBI}.

An amusing variant of the square-ice model is a tiling model
in which the tiles are deformed versions of squares
that the physicist Joshua Burton has dubbed
``gaskets'' and ``baskets'',
depicted in Figure 10.
(To see why the basket deserves its name,
you might want to rotate the page by 45 degrees,
so that the ``handle'' of the basket is pointing up.)
$$
\begin{pspicture}(0.25,-0.25)(1.75,1.25)
\psarc(0,.5){0.707}{-45}{45}
\psarc(2,.5){0.707}{135}{225}
\psarc(1,.5){0.707}{45}{135}
\psarc(1,.5){0.707}{225}{315}
\end{pspicture}
\begin{pspicture}(0,-0.25)(1.75,1.25)
\psarc(0,0.5){0.707}{-45}{45}
\psarc(1,0.5){0.707}{-45}{45}
\psarc(1,0.5){0.707}{225}{315}
\psarc(1,1.5){0.707}{225}{315}
\end{pspicture}
{}
$$
\begin{center} 
Figure 10: A gasket and a basket.
\end{center}
The gasket and basket correspond respectively to
the first and last vertex-types shown in Figure 8;
the other five vertex-types
correspond to tiles obtained by rotating the gasket
by 90 degrees
or by rotating the basket by 90, 180, or 270 degrees.
The directions of the bulges of the four sides of a tile
correspond to the orientations of the four edges
incident with a vertex.
Thus, the seven ASMs of order 3
correspond to the seven distinct ways
of tiling the region shown in Figure 11
(a ``supergasket'' of order 3)
with gaskets and baskets.
$$
\begin{pspicture}(0,0)(2.8,2.8)
\psset{xunit=.7,yunit=.7}
\psset{arrowsize=5pt}
\newcommand{\iceR}{\psarc(0,0){0.4949}{-45}{45}}
\newcommand{\iceL}{\psarc(0,0){0.4949}{135}{225}}
\newcommand{\iceU}{\psarc(0,0){0.4949}{45}{135}}
\newcommand{\iceD}{\psarc(0,0){0.4949}{225}{315}}
\rput(0,1){\iceR}
\rput(0,2){\iceR}
\rput(0,3){\iceR}
\rput(4,1){\iceL}
\rput(4,2){\iceL}
\rput(4,3){\iceL}
\rput(1,3){\iceU}
\rput(2,3){\iceU}
\rput(3,3){\iceU}
\rput(1,1){\iceD}
\rput(2,1){\iceD}
\rput(3,1){\iceD}
\rput(1,3){\iceR}
\rput(1,2){\iceR}
\rput(2,3){\iceR}
\rput(3,2){\iceL}
\rput(2,1){\iceL}
\rput(3,1){\iceL}
\rput(1,2){\iceU}
\rput(1,1){\iceU}
\rput(2,2){\iceU}
\rput(2,2){\iceD}
\rput(3,3){\iceD}
\rput(3,2){\iceD}
\end{pspicture}
\ \ \ \ \ \ \ \ {}
\begin{pspicture}(0,0)(2.8,2.8)
\psset{xunit=.7,yunit=.7}
\psset{arrowsize=5pt}
\newcommand{\iceR}{\psarc(0,0){0.4949}{-45}{45}}
\newcommand{\iceL}{\psarc(0,0){0.4949}{135}{225}}
\newcommand{\iceU}{\psarc(0,0){0.4949}{45}{135}}
\newcommand{\iceD}{\psarc(0,0){0.4949}{225}{315}}
\rput(0,1){\iceR}
\rput(0,2){\iceR}
\rput(0,3){\iceR}
\rput(4,1){\iceL}
\rput(4,2){\iceL}
\rput(4,3){\iceL}
\rput(1,3){\iceU}
\rput(2,3){\iceU}
\rput(3,3){\iceU}
\rput(1,1){\iceD}
\rput(2,1){\iceD}
\rput(3,1){\iceD}
\rput(1,3){\iceR}
\rput(2,3){\iceR}
\rput(1,1){\iceR}
\rput(3,1){\iceL}
\rput(2,2){\iceL}
\rput(3,2){\iceL}
\rput(1,2){\iceU}
\rput(2,2){\iceU}
\rput(2,1){\iceU}
\rput(1,2){\iceD}
\rput(3,3){\iceD}
\rput(3,2){\iceD}
\end{pspicture}
\ \ \ \ \ \ \ \ {}
\begin{pspicture}(0,0)(2.8,2.8)
\psset{xunit=.7,yunit=.7}
\psset{arrowsize=5pt}
\newcommand{\iceR}{\psarc(0,0){0.4949}{-45}{45}}
\newcommand{\iceL}{\psarc(0,0){0.4949}{135}{225}}
\newcommand{\iceU}{\psarc(0,0){0.4949}{45}{135}}
\newcommand{\iceD}{\psarc(0,0){0.4949}{225}{315}}
\rput(0,1){\iceR}
\rput(0,2){\iceR}
\rput(0,3){\iceR}
\rput(4,1){\iceL}
\rput(4,2){\iceL}
\rput(4,3){\iceL}
\rput(1,3){\iceU}
\rput(2,3){\iceU}
\rput(3,3){\iceU}
\rput(1,1){\iceD}
\rput(2,1){\iceD}
\rput(3,1){\iceD}
\rput(1,2){\iceR}
\rput(1,3){\iceR}
\rput(2,2){\iceR}
\rput(3,1){\iceL}
\rput(3,3){\iceL}
\rput(2,1){\iceL}
\rput(1,1){\iceU}
\rput(1,2){\iceU}
\rput(3,1){\iceU}
\rput(3,2){\iceU}
\rput(2,3){\iceD}
\rput(2,2){\iceD}
\end{pspicture}
\ \ \ \ \ \ \ \ {}
\begin{pspicture}(0,0)(2.8,2.8)
\psset{xunit=.7,yunit=.7}
\psset{arrowsize=5pt}
\newcommand{\iceR}{\psarc(0,0){0.4949}{-45}{45}}
\newcommand{\iceL}{\psarc(0,0){0.4949}{135}{225}}
\newcommand{\iceU}{\psarc(0,0){0.4949}{45}{135}}
\newcommand{\iceD}{\psarc(0,0){0.4949}{225}{315}}
\rput(0,1){\iceR}
\rput(0,2){\iceR}
\rput(0,3){\iceR}
\rput(4,1){\iceL}
\rput(4,2){\iceL}
\rput(4,3){\iceL}
\rput(1,3){\iceU}
\rput(2,3){\iceU}
\rput(3,3){\iceU}
\rput(1,1){\iceD}
\rput(2,1){\iceD}
\rput(3,1){\iceD}
\rput(1,3){\iceR}
\rput(2,2){\iceR}
\rput(1,1){\iceR}
\rput(3,3){\iceL}
\rput(2,2){\iceL}
\rput(3,1){\iceL}
\rput(1,2){\iceU}
\rput(2,1){\iceU}
\rput(3,2){\iceU}
\rput(1,2){\iceD}
\rput(2,3){\iceD}
\rput(3,2){\iceD}
\end{pspicture}
$$
$$
\begin{pspicture}(0,0)(2.8,2.8)
\psset{xunit=.7,yunit=.7}
\psset{arrowsize=5pt}
\newcommand{\iceR}{\psarc(0,0){0.4949}{-45}{45}}
\newcommand{\iceL}{\psarc(0,0){0.4949}{135}{225}}
\newcommand{\iceU}{\psarc(0,0){0.4949}{45}{135}}
\newcommand{\iceD}{\psarc(0,0){0.4949}{225}{315}}
\rput(0,1){\iceR}
\rput(0,2){\iceR}
\rput(0,3){\iceR}
\rput(4,1){\iceL}
\rput(4,2){\iceL}
\rput(4,3){\iceL}
\rput(1,3){\iceU}
\rput(2,3){\iceU}
\rput(3,3){\iceU}
\rput(1,1){\iceD}
\rput(2,1){\iceD}
\rput(3,1){\iceD}
\rput(1,1){\iceR}
\rput(2,1){\iceR}
\rput(1,3){\iceR}
\rput(3,3){\iceL}
\rput(2,2){\iceL}
\rput(3,2){\iceL}
\rput(1,2){\iceU}
\rput(3,2){\iceU}
\rput(3,1){\iceU}
\rput(1,2){\iceD}
\rput(2,2){\iceD}
\rput(2,3){\iceD}
\end{pspicture}
\ \ \ \ \ \ \ \ {}
\begin{pspicture}(0,0)(2.8,2.8)
\psset{xunit=.7,yunit=.7}
\psset{arrowsize=5pt}
\newcommand{\iceR}{\psarc(0,0){0.4949}{-45}{45}}
\newcommand{\iceL}{\psarc(0,0){0.4949}{135}{225}}
\newcommand{\iceU}{\psarc(0,0){0.4949}{45}{135}}
\newcommand{\iceD}{\psarc(0,0){0.4949}{225}{315}}
\rput(0,1){\iceR}
\rput(0,2){\iceR}
\rput(0,3){\iceR}
\rput(4,1){\iceL}
\rput(4,2){\iceL}
\rput(4,3){\iceL}
\rput(1,3){\iceU}
\rput(2,3){\iceU}
\rput(3,3){\iceU}
\rput(1,1){\iceD}
\rput(2,1){\iceD}
\rput(3,1){\iceD}
\rput(1,1){\iceR}
\rput(1,2){\iceR}
\rput(2,2){\iceR}
\rput(3,3){\iceL}
\rput(2,3){\iceL}
\rput(3,1){\iceL}
\rput(2,2){\iceU}
\rput(2,1){\iceU}
\rput(3,2){\iceU}
\rput(1,2){\iceD}
\rput(1,3){\iceD}
\rput(3,2){\iceD}
\end{pspicture}
\ \ \ \ \ \ \ \ {}
\begin{pspicture}(0,0)(2.8,2.8)
\psset{xunit=.7,yunit=.7}
\psset{arrowsize=5pt}
\newcommand{\iceR}{\psarc(0,0){0.4949}{-45}{45}}
\newcommand{\iceL}{\psarc(0,0){0.4949}{135}{225}}
\newcommand{\iceU}{\psarc(0,0){0.4949}{45}{135}}
\newcommand{\iceD}{\psarc(0,0){0.4949}{225}{315}}
\rput(0,1){\iceR}
\rput(0,2){\iceR}
\rput(0,3){\iceR}
\rput(4,1){\iceL}
\rput(4,2){\iceL}
\rput(4,3){\iceL}
\rput(1,3){\iceU}
\rput(2,3){\iceU}
\rput(3,3){\iceU}
\rput(1,1){\iceD}
\rput(2,1){\iceD}
\rput(3,1){\iceD}
\rput(1,1){\iceR}
\rput(1,2){\iceR}
\rput(2,1){\iceR}
\rput(3,3){\iceL}
\rput(2,3){\iceL}
\rput(3,2){\iceL}
\rput(2,2){\iceU}
\rput(3,1){\iceU}
\rput(3,2){\iceU}
\rput(1,2){\iceD}
\rput(1,3){\iceD}
\rput(2,2){\iceD}
\end{pspicture}
$$
\begin{center} 
Figure 11: The seven tilings of an order-3 supergasket 
with gaskets and baskets.
\end{center}

\noindent
For another fanciful embodiment of ASMs as tilings,
see the cover of~\cite{Br}. 

\section{Symmetric ASMs and partial ASMs}
\label{sec:sym}
Some ASMs are more symmetrical than others.
More precisely,
the eight-element dihedral group $D_4$
acts on ASMs,
and for every subgroup $G$ of $D_4$
there are ASMs that are invariant under the action
of every element of $G$.
In~\cite{Ro2},
Robbins gave some conjectures for
the number of ASMs of order $n$
that are invariant under particular groups $G$;
for most (but not all) of the subgroups $G$ of $D_4$,
numerical evidence suggested specific product-formula.
Since then,
Kuperberg~\cite{Ku2}
has proved some of these, but others remain conjectural.

At the same time,
one may also look at halves
(or even quarters or eighths) of ASMs ---
the fundamental regions under the action
of the aforementioned groups $G$ ---
and look at them in their own right,
asking, How many partial ASMs are there
if one limits attention to such a region?
There are some interesting phenomena here.

For instance, for $c_1$,$c_2$,$c_3$
each equal to $+1$ or $-1$,
define $N(c_1,c_2,c_3)$ as the number of 
4-by-7 partial height-function matrices
of the form shown in Figure 12.
$$
\begin{array}{ccccccc}
0 & 1 & 2 & 3 & 4 & 5 & 6 \\
1 & ? & ? & ? & ? & ? & 5 \\
2 & ? & ? & ? & ? & ? & 4 \\
3 & 3+c_1 & 3 & 3+c_2 & 3 & 3+c_3 & 3
\end{array}
$$
\begin{center}
Figure 12: Half of a height-function matrix of order 6.
\end{center}
Not surprisingly, the eight values of $N(c_1,c_2,c_3)$
as $c_1,c_2,c_3$ vary
are not all equal to one another.
But it is surprising that the four numbers
$$N(1,1,1),$$
$$(N(1,1,-1)+N(1,-1,1)+N(-1,1,1))/3,$$
$$(N(1,-1,-1)+N(-1,1,-1)+N(-1,-1,1))/3,$$
and 
$$N(-1,-1,-1)$$ 
{\it are\/} all equal.
More generally,
consider $(n+1)$-by-$(2n+1)$ partial height-function matrices
of the following form:
$$
\begin{array}{ccccccccc}
0 & 1 & 2 & 3 & 4 & 5 & ... & 2n-1 & 2n \\
1 & ? & ? & ? & ? & ? & ... & ? & 2n-1 \\
2 & ? & ? & ? & ? & ? & ... & ? & 2n-2 \\
\vdots & \vdots & \vdots & \vdots & \vdots & \vdots & \ddots & \vdots & \vdots \\
n-1 & ? & ? & ? & ? & ? & ... & ? & n+1 \\
n & n+c_1 & n & n+c_2 & n & n+c_3 & ... & n+c_{n} & n
\end{array} .
$$
\begin{center}
Figure 13: Half of a height-function matrix of order $2n$.
\end{center}
Here each $c_i$ ($1 \leq i \leq n$) is either $+1$ or $-1$.
Let $N(c_1,c_2,\dots,c_n)$ be the number of such
partial height-function matrices.
Then one finds empirically that
for every $k$ in $-n,-(n-2),-(n-4),\dots,n-4,n-2,n$,
the average of $N(c_1,c_2,\dots,c_n)$ over all vectors $(c_1,\dots,c_n)$
satisfying $c_1+\dots+c_n=k$
depends only on $n$, not on $k$.
Kuperberg has found an algebraic proof of this 
using the Tsuchiya determinant formula~\cite{T} 
invented for the study of the square ice model,
but there ought to be a purely combinatorial proof
of this simple relation.

\section{Weighted enumeration}
\label{sec:weight}
There are some interesting results in the literature
on weighted enumeration of ASMs.
Here one assigns to each ASM of order $n$ some weight,
and tries to compute the sum of the weights of all the ASMs of order $n$.
A priori it might be unclear why this would be interesting,
but with certain weighting schemes
one gets beautiful (and mysterious) formulas,
which are their own justification.

For instance, following~\cite{MRR}, 
one can assign weight $x^k$ to every ASM
that contains exactly $k$ entries equal to $-1$.
What is the sum of the weights of the ASMs of order $n$?
When $x=1$, this is nothing other than ordinary enumeration of ASMs.
When $x=2$, there is a very nice answer~\cite{MRR}:
the sum of the weights is exactly $2^{n(n-1)/2}$.
When $x=3$, the answer is more complicated,
but it is roughly similar in type to the formula for the case $x=1$,
and roughly similar in difficulty;
the ``3-enumeration'' formula
was first conjectured by Mills, Robbins and Rumsey~\cite{MRR}
and was eventually proved by Kuperberg~\cite{Ku2}
(with corrections provided by Robin Chapman).
No other positive integer $x$ seems to give nice answers.
One can also assign weight $x^k$ to every ASM
that contains exactly $k$ entries equal to $+1$,
but this is essentially the same weighting scheme,
since in any ASM of order $n$,
the number of $+1$'s is always
$n$ plus the number of $-1$'s.

More interestingly, one can also use a hybrid weighting scheme
in which the exponent of $x$ is equal to the number of entries
$a_{i,j}$ such that \emph{either} $i+j$ is even and $a_{i,j}=-1$
\emph{or} $i+j$ is odd and $a_{i,j}=+1$.
When $x=2$, this too leads to an interesting result:
the sum of the weights is always a power of two
times a power of five!~\cite{Y}

One can come up with many open problems by combining
the ideas of this section and the previous section.
Here is one example:
Each way of filling in Figure 13
(with the $c_i$'s now permitted to vary freely)
gives rise to a ``half-ASM''.
For instance, the filling
$$
\begin{array}{ccccccc}
0 & 1 & 2 & 3 & 4 & 5 & 6 \\
1 & 2 & 3 & 4 & 5 & 4 & 5 \\
2 & 3 & 4 & 3 & 4 & 3 & 4 \\
3 & 2 & 3 & 2 & 3 & 4 & 3
\end{array}
$$
of Figure 12 gives rise to the half-ASM
$$
\begin{array}{rrrrrr}
\ 0 & \ \ 0 & \ 0 & \ \ 0 &  +1 & \ 0 \\
\ 0 & \ \ 0 &  +1 & \ \ 0 & \ 0 & \ 0 \\
 +1 & \ \ 0 & \ 0 & \ \ 0 &  -1 &  +1 
\end{array}
$$
which has a single $-1$.
If we assign each rectangular array that arises in this way
a weight equal to 2 to the power of the number of $-1$'s, 
we get $2^{n^2}$ (Robin Chapman and Theresia Eisenk\"olbl 
have independently found nice proofs of this).
On the other hand, suppose we now permit
every entry in the bottom row of Figure 13
to vary freely
(aside from the $n$ on the left and the $n$ on the right).
When we 2-enumerate half-ASMs of this sort,
as a function of $n$,
we get the following sequence of numbers:
$2$, $20 = 2^2 \cdot 5$,
$896 = 2^7 \cdot 7$,
$177408 = 2^8 \cdot 3^2 \cdot 7 \cdot 11$,
$154632192 = 2^{15} \cdot 3 \cdot 11^2 \cdot 13$,
$592344383488 = 2^{17} \cdot 11^2 \cdot 13^3 \cdot 17$, $\dots$.
Clearly the absence of larger prime factors
indicates that there is some nice product formula
governing these numbers.
Can someone find the right conjecture?
Can someone prove it?
(For more data of this kind,
see: 
\begin{center}
{\tt http://www.math.wisc.edu/$\sim$propp/half-asm}
\end{center}
Late-breaking news: Theresia Eisenk\"olbl~\cite{E} has made
progress with the data-set
and proved a number of theorems about half-ASMs.)

\section{Full packings of loops}
\label{sec:loops}
Given an ice state of order $n$,
we can form a subgraph of the underlying tic-tac-toe graph
by selecting precisely those edges
that are oriented so as to point from an odd vertex to an even vertex,
where we have assigned parities to vertices
so that each odd vertex has only even vertices as neighbors
and vice versa.
Then one gets a subgraph of the tic-tac-toe graph
such that each of the $n^2$ internal vertices
lies on exactly 2 of the selected edges,
and the $4n$ external vertices,
taken in cyclical order,
alternate between lying on a selected edge
and not lying on a selected edge.
Moreover, every such subgraph
arises from an ice state in this way.

Let us say that the leftmost vertex
in the top row of external vertices
is even.
Figure 14 shows the seven subgraphs that result
from applying the transformation
to the seven ice states of order 3.

$$
\begin{pspicture}(0,0)(3.6,3.6)
\psset{xunit=.6,yunit=.6}
\multips(2,1)(1,0){3}{\odisk{.1}}
\multips(1,2)(1,0){5}{\odisk{.1}}
\multips(1,3)(1,0){5}{\odisk{.1}}
\multips(1,4)(1,0){5}{\odisk{.1}}
\multips(2,5)(1,0){3}{\odisk{.1}}
\psline(2,4)(2,5)
\psline(4,4)(4,5)
\psline(2,1)(2,2)
\psline(4,1)(4,2)
\psline(1,3)(2,3)
\psline(4,3)(5,3)
\psline(2,4)(3,4)(3,2)(4,2)
\psline(4,4)(4,3)
\psline(2,2)(2,3)
\end{pspicture}
\ {}
\begin{pspicture}(0,0)(4.2,4.2)
\psset{xunit=.6,yunit=.6}
\multips(2,1)(1,0){3}{\odisk{.1}}
\multips(1,2)(1,0){5}{\odisk{.1}}
\multips(1,3)(1,0){5}{\odisk{.1}}
\multips(1,4)(1,0){5}{\odisk{.1}}
\multips(2,5)(1,0){3}{\odisk{.1}}
\psline(2,4)(2,5)
\psline(4,4)(4,5)
\psline(2,1)(2,2)
\psline(4,1)(4,2)
\psline(1,3)(2,3)
\psline(4,3)(5,3)
\psline(2,4)(3,4)(3,3)(2,3)
\psline(4,4)(4,3)
\psline(2,2)(4,2)
\end{pspicture}
\ {}
\begin{pspicture}(0,0)(4.2,4.2)
\psset{xunit=.6,yunit=.6}
\multips(2,1)(1,0){3}{\odisk{.1}}
\multips(1,2)(1,0){5}{\odisk{.1}}
\multips(1,3)(1,0){5}{\odisk{.1}}
\multips(1,4)(1,0){5}{\odisk{.1}}
\multips(2,5)(1,0){3}{\odisk{.1}}
\psline(2,4)(2,5)
\psline(4,4)(4,5)
\psline(2,1)(2,2)
\psline(4,1)(4,2)
\psline(1,3)(2,3)
\psline(4,3)(5,3)
\psline(2,4)(4,4)
\psline(2,2)(2,3)
\psline(4,2)(3,2)(3,3)(4,3)
\end{pspicture}
\ {}
\begin{pspicture}(0,0)(4.2,4.2)
\psset{xunit=.6,yunit=.6}
\multips(2,1)(1,0){3}{\odisk{.1}}
\multips(1,2)(1,0){5}{\odisk{.1}}
\multips(1,3)(1,0){5}{\odisk{.1}}
\multips(1,4)(1,0){5}{\odisk{.1}}
\multips(2,5)(1,0){3}{\odisk{.1}}
\psline(2,4)(2,5)
\psline(4,4)(4,5)
\psline(2,1)(2,2)
\psline(4,1)(4,2)
\psline(1,3)(2,3)
\psline(4,3)(5,3)
\psline(2,4)(4,4)
\psline(2,3)(4,3)
\psline(2,2)(4,2)
\end{pspicture}
$$
$$
\begin{pspicture}(0,0)(4.2,4.2)
\psset{xunit=.6,yunit=.6}
\multips(2,1)(1,0){3}{\odisk{.1}}
\multips(1,2)(1,0){5}{\odisk{.1}}
\multips(1,3)(1,0){5}{\odisk{.1}}
\multips(1,4)(1,0){5}{\odisk{.1}}
\multips(2,5)(1,0){3}{\odisk{.1}}
\psline(2,4)(2,5)
\psline(4,4)(4,5)
\psline(2,1)(2,2)
\psline(4,1)(4,2)
\psline(1,3)(2,3)
\psline(4,3)(5,3)
\psline(2,4)(4,4)
\psline(2,3)(3,3)(3,2)(2,2)
\psline(4,2)(4,3)
\end{pspicture}
\ {}
\begin{pspicture}(0,0)(4.2,4.2)
\psset{xunit=.6,yunit=.6}
\multips(2,1)(1,0){3}{\odisk{.1}}
\multips(1,2)(1,0){5}{\odisk{.1}}
\multips(1,3)(1,0){5}{\odisk{.1}}
\multips(1,4)(1,0){5}{\odisk{.1}}
\multips(2,5)(1,0){3}{\odisk{.1}}
\psline(2,4)(2,5)
\psline(4,4)(4,5)
\psline(2,1)(2,2)
\psline(4,1)(4,2)
\psline(1,3)(2,3)
\psline(4,3)(5,3)
\psline(2,4)(2,3)
\psline(2,2)(4,2)
\psline(4,4)(3,4)(3,3)(4,3)
\end{pspicture}
\ {}
\begin{pspicture}(0,0)(4.2,4.2)
\psset{xunit=.6,yunit=.6}
\multips(2,1)(1,0){3}{\odisk{.1}}
\multips(1,2)(1,0){5}{\odisk{.1}}
\multips(1,3)(1,0){5}{\odisk{.1}}
\multips(1,4)(1,0){5}{\odisk{.1}}
\multips(2,5)(1,0){3}{\odisk{.1}}
\psline(2,4)(2,5)
\psline(4,4)(4,5)
\psline(2,1)(2,2)
\psline(4,1)(4,2)
\psline(1,3)(2,3)
\psline(4,3)(5,3)
\psline(2,4)(2,3)
\psline(4,3)(4,2)
\psline(4,4)(3,4)(3,2)(2,2)
\end{pspicture}
$$
\begin{center}
Figure 14: The seven FPL states of order 3.
\end{center}

Leaving aside the behavior at the boundary,
these are states of what physicists call
the fully packed loop (FPL) model
on the square grid (see e.g.~\cite{BBNY}).
I sometimes prefer to call such states
``near 2-factors''
since nearly all of the vertices in these subgraphs
have degree 2;
only the external vertices have smaller degree
(specifically, they alternate between having degree 0
and having degree 1).

If one starts from an external vertex,
there is a unique path that one can follow
using edges in the subgraph;
this path must eventually lead to one of the other external vertices.
In addition to these paths (``open loops''),
the edges of the subgraph can also form closed loops 
(see Figure 15, for example).

Note that these loops (both open and closed)
cannot cross one another.
In particular, the open loops
must join up the $2n$ external edges
in some non-crossing fashion.
If one numbers the vertices of degree 1
in cyclic order from 1 to $2n$,
the FPL state yields a pairing
of odd-indexed external vertices
with even-indexed external vertices.
For instance,
the FPL state shown in Figure 15
links 1 with 12,
2 with 11,
3 with 4,
5 with 6,
7 with 8,
and 9 with 10.
$$
\begin{pspicture}(.4,.4)(5.0,5.0)
\psset{xunit=.6,yunit=.6}
\multips(2,1)(1,0){6}{\odisk{.1}}
\multips(1,2)(1,0){8}{\odisk{.1}}
\multips(1,3)(1,0){8}{\odisk{.1}}
\multips(1,4)(1,0){8}{\odisk{.1}}
\multips(1,5)(1,0){8}{\odisk{.1}}
\multips(1,6)(1,0){8}{\odisk{.1}}
\multips(1,7)(1,0){8}{\odisk{.1}}
\multips(2,8)(1,0){6}{\odisk{.1}}
\psline(1,6)(2,6)(2,8)
\psline(1,4)(2,4)(2,5)(3,5)(3,7)(4,7)(4,8)
\psline(1,2)(2,2)(2,3)(3,3)(3,4)(4,4)(4,2)(3,2)(3,1)
\psline(6,8)(6,7)(5,7)(5,6)(4,6)(4,5)(6,5)(6,6)(7,6)(7,7)(8,7)
\psline(8,5)(7,5)(7,3)(8,3)
\psline(5,1)(5,2)(7,2)(7,1)
\psline(5,3)(6,3)(6,4)(5,4)(5,3)
\rput(2,8.5){1}
\rput(4,8.5){2}
\rput(6,8.5){3}
\rput(8.5,7){4}
\rput(8.5,5){5}
\rput(8.5,3){6}
\rput(7,0.5){7}
\rput(5,0.5){8}
\rput(3,0.5){9}
\rput(0.5,2){10}
\rput(0.5,4){11}
\rput(0.5,6){12}
\end{pspicture}
$$
\begin{center}
Figure 15: A fully packed loop state of order 6.
\end{center}

It has been conjectured, on the strength of numerical evidence,
that the number of ASMs of order $n$
in which the open paths link
1 with 2, 3 with 4, \dots, and $2n-1$ with $2n$
is exactly the total number of ASMs of order $n-1$.

We do not have a proof of this, but curiously, we have a proof
of something else:
that the number of FPL states of order $n$
in which the open paths link
1 with 2, 3 with 4, \dots, and $2n-1$ with $2n$
is equal to the number of FPL states of order $n$
in which the open paths link
1 with $2n$, 2 with 3, 4 with 5, \dots, and $2n-2$ with $2n-1$.
Note that when $n$ is divisible by 4,
the geometries of the two linking-patterns is different,
with respect to the tic-tac-toe graph.
Yet the number of FPL states is the same.

This is a special case of a far more general fact
proved by Wieland~\cite{W}.
For any two non-crossing pairings $\pi$ and $\pi'$
of the numbers 1 through $2n$
(viewed as equally spaced points on a circle),
if $\pi$ and $\pi'$ are conjugate
via a rotation or reflection,
then (if we now treat the numbers 1 through $2n$
as the labels of vertices of degree 1
in the tic-tac-toe graph of order $n$)
the number of FPL states with linking-pattern $\pi$ 
equals the number of FPL states with linking-pattern $\pi'$.
It is as if the tic-tac-toe graph,
in some mystical sense,
had an automorphism sending 1 to 2, 2 to 3, etc.

One might also ask for the number of FPL states of order $n$
in which 1 is linked with 2
(ignoring all the other linking going on).
Here too we have a conjectural answer,
due to David Wilson:
the number of such FPL states 
is just the total number of FPL states of order $n$
multiplied by
$$\frac{3}{2} \frac{n^2+1}{4n^2-1}.$$
This would imply, in particular,
that as $n$ goes to infinity,
the probability that a randomly chosen FPL state
links 1 with 2 
is asymptotically $3/8$.

For a very recent discussion of the FPL model, see~\cite{RS2}.
In this article, Razumov and Stroganov point out that
the FPL model is closely related to a seemingly quite
different lattice model.

It is worth pointing out that these conjectures
are truly native to the FPL incarnation of ASMs;
it is hard to see how they could have arisen from
one of the other models.
So, even though the transformation
between ASMs and FPL states is fairly shallow mathematically,
some interesting questions can arise from it
that might not otherwise have been noticed.
Likewise, the passage between square-ice states
and ASMs is not conceptually deep,
but it made possible the shortest known solution of the ASM problem
by putting the problem into a form
where known methods from the physics literature could be applied.
Hence a good subtitle for this paper might have been
``The non-trivial power of trivial transformations''.

\section{Numerology}
\label{sec:numerology}
As far as I know,
the first manifestation of the sequence
1,2,7,42,429,7436,\dots
occurred in connection with combinatorial objects
called descending plane partitions or DPPs~\cite{An1}.
Another manifestation was totally symmetric
self-complementary plane-partitions (TSSCPPs)~\cite{An2}.
Andrews discovered the proofs of both these formulas.
Indeed, it was Andrews' proof of the 1,2,7,42,\dots formula
for TSSCPPs
that galvanized Zeilberger into tackling the ASM conjecture.
Zeilberger showed that ASMs are equinumerous with TSSCPPs.
However, his proof was not bijective,
and to this day nobody knows of a good bijection
between ASMs of order $n$
and TSSCPPs of order $n$.

Another context in which these numbers arise
is the study of the XXZ model 
in statistical mechanics~\cite{RS1}~\cite{BGN}.

A different context in which numbers related to ASMs have occurred
is certain ``number walls'' investigated by Somos~\cite{So}.
(Number walls are arrays of Hankel determinants,
arranged so as to faciliate calculations of
successively larger ones;
see~\cite{CG} for details.)
In Somos' examples
it is not the sequence 1,2,7,42,\dots that crops up,
but sequences that enumerate various sorts of symmetric ASMs.
Xin~\cite{X} has found a Hankel determinant theorem
that involves the sequence $1,2,7,42,\dots$ itself:
he has shown that for all $n$, the determinant of the $n$-by-$n$ 
matrix whose $i,j$th entry is equal to the coefficient of
$x^{i+j-2}$ in the Taylor expansion of
the generalized Catalan generating function
$\frac{1-(1-9x)^{1/3}}{3x}$~\cite{La}
is equal to
$$3^{n \choose 2}$$
times the number of $n$-by-$n$ alternating-sign matrices.

It would be desirable to have some sort of understanding
of why the number of ASMs of order $n$
(with or without symmetry-constraints)
turns up in these seemingly disparate situations.

\section{Large random ASMs}
\label{sec:large}
Another sort of phenomenon associated with ASMs of order $n$
is their typical ``shape'' when $n$ is large.
I remarked above that the ASMs of order $n$
form a distributive lattice;
consequently, the method of ``coupling from the past''
can be applied~\cite{Pr}
for the purpose of generating 
a random ASM of fixed order
governed by the uniform distribution.
Figure 16 shows a random ASM of order 40,
represented as a gaskets-and-baskets tiling
(where an attempt has been made to give
each of the six tile-types its own distinctive shading).

\begin{figure}
\begin{center}
\leavevmode
\epsfbox[100 100 350 500]{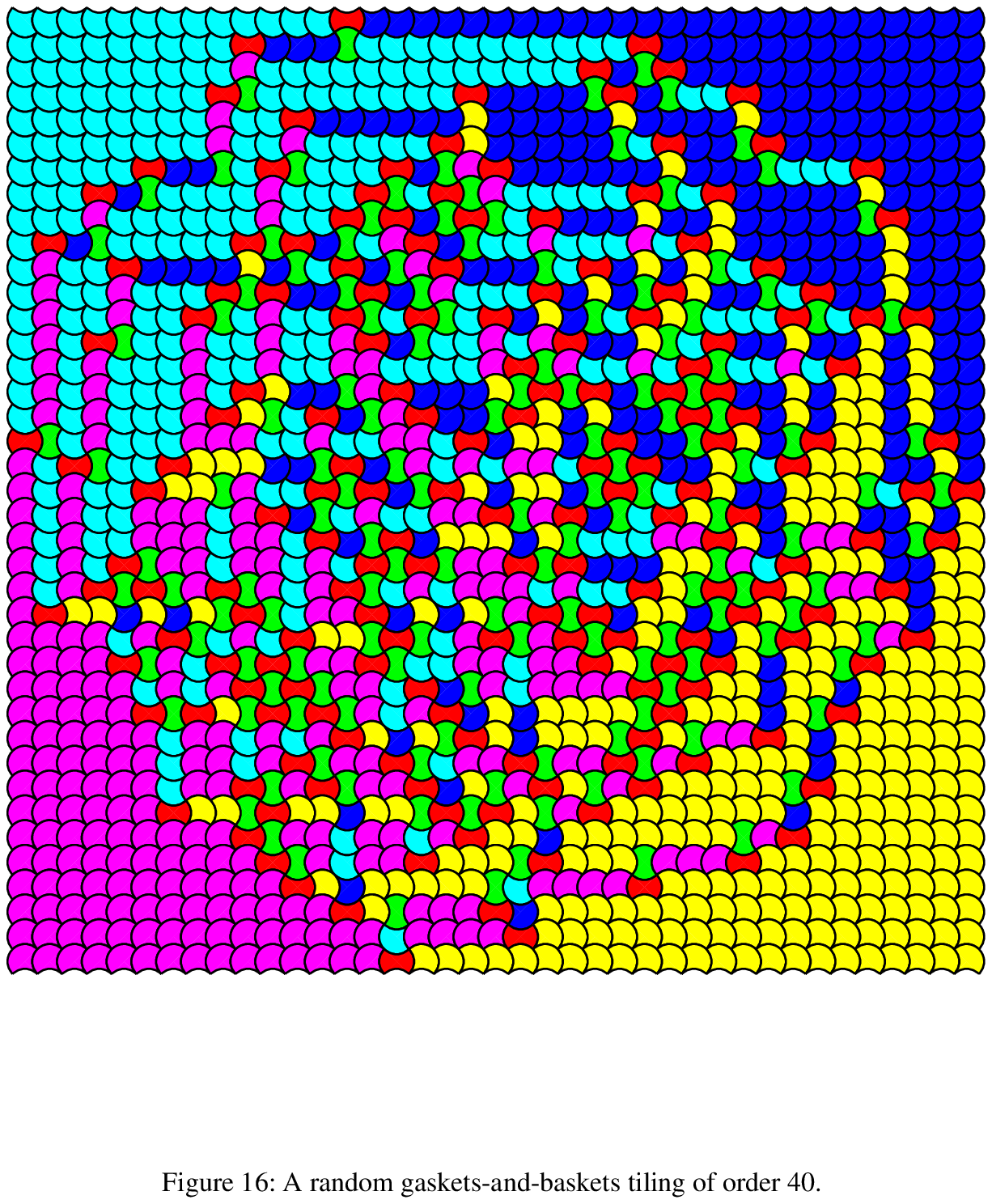}
\end{center}
\end{figure}

Note that the gaskets (which correspond to
the non-zero entries of the ASM)
stay away from the corners.
Computer experiments strongly indicate
that this is typical behavior:
the probability of finding a non-zero entry
close to one of the corners
appears to be quite small.
Another way of expressing this is
in terms of the entries $h_{i,j}$ of the height-function matrix.
Say that a location $(i,j)$ in a particular height-function matrix
of order $n$
is \emph{frozen}
if the height there is equal to either
the maximum possible height 
that any height-function matrix of order $n$ can exhibit at that location
or the minimum possible height 
that any height-function matrix of order $n$ can exhibit at that location.
Then the claim is that a significant portion of the height-function
matrix, concentrated in the four corners,
tends to be frozen.

This is analogous to phenomena that have been observed
for other sorts of combinatorial models.
Indeed, if one adopt a non-uniform distribution
on the set of ASMs of order $n$,
where the probability associated with an ASM
containing exactly $k$ entries equal to $-1$
is proportional to $2^k$,
then it is rigorously known
that the frozen region
tend in probability to a perfect circular disk~\cite{JPS}~\cite{CEP}.
However, it is not known that this holds
when the uniform distribution on ASMs is used.

\section{Back to Dodgson}
\label{sec:dodgson}
I conclude this article by coming full circle
and returning the context in which
ASMs first came to light:
the study of Dodgson's condensation algorithm and its variants.
I will not discuss Dodgson's algorithm per se,
but rather a variation of it invented by Robbins and Rumsey~\cite{RR}.
This modified form of Dodgson condensation
is an algebraic recurrence relation
$$f_{i,j,k+1}=(f_{i+1,j,k}f_{i-1,j,k}+f_{i,j+1,k}f_{i,j-1,k})/f_{i,j,k-1}$$
satisfied by certain functions $f:\Z^3 \rightarrow \R$.
(Note that this equation can be written slightly more symmetrically as
$$f_{i,j,k+1}f_{i,j,k-1}-f_{i+1,j,k}f_{i-1,j,k}-f_{i,j+1,k}f_{i,j-1,k}=0.$$
Physicists and researchers in the field of integrable systems
call this a discrete Hirota equation, 
and have developed a great deal of theory associated with it;
however, the observations I make here
seem to be currently unknown 
outside of a small circle of algebraic combinatorialists.)
If we let $f_{i,j,-1}=x_{i,j}$
and $f_{i,j,0}=y_{i,j}$
for formal indeterminates $x_{i,j},y_{i,j}$
(with $(i,j)$ ranging over $\Z^2$),
then the recurrence relation lets us express
all the $f_{i,j,k}$'s
in terms of the $x_{i,j}$'s and $y_{i,j}$'s,
at least formally.
A priori, one expects each $f_{i,j,k}$ to be
a rational function of the $x$ and $y$ variables;
the surprise (the first of several surprises, in fact)
is that these rational functions are actually Laurent polynomials
(that is, they are polynomials functions of
the $x$ and $y$ variables along with their reciprocals).

This observation seems to have first been made by Mills, Robbins and Rumsey,
who actually considered a more general recurrence
$$f_{i,j,k+1}=(f_{i+1,j,k}f_{i-1,j,k} + 
\lambda f_{i,j+1,k}f_{i,j-1,k})/f_{i,j,k-1}.$$
The case $\lambda=-1$ corresponds to the original Dodgson algorithm,
but Mills et al.\ noticed that the same surprising cancellation occurs
for more general values of $\lambda$,
including the especially nice case $\lambda=1$,
and that one always obtains Laurent polynomials.

It was from studying these Laurent polynomials
that Mills, Robbins and Rumsey were led 
to discover alternating-sign matrices.
Every term in one of these Laurent polynomials
has a coefficient equal to 1 (that is the second surprise),
and is a product of powers of a finite number
of $x$ and $y$ variables.
The third surprise is that all the exponents of the variables
are $+1$, $-1$, and $0$.
The fourth and final surprise is that
these patterns of exponents encode ASMs.
More specifically, the exponents of the $x$ variables
(after a global sign flip)
encode one ASM,
and the exponents of the $y$ variables
encode another.
These two ASMs satisfy a combinatorial relationship 
that the researchers dubbed ``compatibility''.
They showed that the number of compatible pairs of ASMs
is exactly $2^{n(n+1)/2}$.

As it happens, this formula is not hard to verify
as a consequence of the other claims I have made.
If indeed all the coefficients in the Laurent polynomial equal 1,
then one can count the terms
(and thereby count the compatible pairs of ASMs)
just by setting all the $x$ and $y$ variables equal to 1.
But in this case, $f_{i,j,k}$ depends only on $k$ (call it $F_k$),
and the three-dimensional recurrence boils down 
to the one-dimensional recurrence
$$F_{k+1}=(F_k F_k + F_k F_k)/F_{k-1}$$
with initial conditions $F_0=F_1=1$,
which is readily solved.

The terms of these Laurent polynomials 
were originally understood in terms of
compatible pairs of ASMs.
A few years after this work was done,
it turned out that compatible pairs of ASMs
admit a much more geometrical representation,
namely, as tilings of regions called Aztec diamonds
by means of tiles called dominos (1-by-2 and 2-by-1 rectangles).
See~\cite{EKLP} for more details.

I will close by pointing out that a kindred recurrence relation
cries out to be studied, namely
$$f_{i,j,k} =(f_{i-1,j,k}f_{i,j-1,k-1}+f_{i,j-1,k}f_{i-1,j,k-1}
+f_{i,j,k-1}f_{i-1,j-1,k})/f_{i-1,j-1,k-1},$$
with initial conditions
$$f_{i,j,k}=
\left\{ \begin{array}{ll}
x_{i,j,k} & \mbox{if $i+j+k=-1$}, \\
y_{i,j,k} & \mbox{if $i+j+k= 0$}, \\
z_{i,j,k} & \mbox{if $i+j+k=+1$}. 
\end{array}
\right.
$$
Here (just as in the Mills, Robbins, and Rumsey recurrence)
one finds empirically
that each value of $f_{i,j,k}$
is expressible as a Laurent polynomial
in the $x$, $y$, and $z$ variables
(in fact, shortly before this article went to press,
this Laurent property was proved by Fomin and Zelevinsky\cite{FZ});
here too one finds empirically
that each coefficient in these Laurent polynomials
is equal to 1;
and here too one finds that the exponents of the $x$, $y$ and $z$ variables
that occur in the Laurent monomials
are universally bounded
(in this case between $-1$ and $+4$
rather than between $-1$ and $+1$).
If all this is true,
then the exponent-patterns that arise
are some sort of analogue of
compatible pairs of ASMs,
and moreover, we know exactly how many there are:
$3^{\lfloor n^2/4 \rfloor}$.
(This comes from reducing the original three-dimensional recurrence
to a one-dimensional recurrence, as we did before.)
So, assuming that our empirical observations are not leading us astray,
there is some new kind of combinatorial gadget
that governs these Laurent polynomials
(or vice versa!),
and we know exactly how many gadgets of order $n$ there are. 
And it is easy to generate these Laurent polynomials 
(and with them the gadgets)
using {\sc Maple},
e.g.\ with the following short program:
\begin{verbatim}
f := proc (i,j,k)
        if   (i+j+k < 3) then x(i,j,k) else
        simplify(
        ( f(i-1,j,k)*f(i,j-1,k-1)+
          f(i,j-1,k)*f(i-1,j,k-1)+
          f(i,j,k-1)*f(i-1,j-1,k) )
        /f(i-1,j-1,k-1));
        fi; end;
\end{verbatim}
Nonetheless, we do not know what these gadgets are, combinatorially!
They are analogous to pairs of compatible ASMs,
which in turn are equivalent to domino tilings of Aztec diamonds,
so one hopes that the gadgets have some intuitive geometric meaning.
For more information about the properties of these gadgets,
see {\tt http://www.math.wisc.edu/$\sim$propp/cube-recur}.

More late-breaking news: Harvard undergraduates Gabriel Carroll
and David Speyer solved this problem in 2002, as participants
in the Research Experiences in Algebraic Combinatorics at
Harvard program.  More details will at some point become
available through links at the REACH web-site:
\begin{center}
{\tt http://www.math.harvard.edu/$\sim$propp/reach/}.
\end{center}


\end{document}